\documentclass{amsart}
\usepackage[latin1]{inputenc}
\usepackage[dvips]{graphicx}
\usepackage{amsfonts,amsmath,amsthm,amssymb,amscd}
\usepackage{enumerate}

\newcommand{\C}{\mathbb{C}}
\newcommand{\R}{\mathbb{R}}
\newcommand{\Z}{\mathbb{Z}}
\newcommand{\N}{\mathbb{N}}

\newcommand{\interior}[1]{\stackrel{\circ}{#1}}

\newcommand{\fl}{\mathit{fl}}
\theoremstyle{plain}
\newtheorem{Lem}{Lemma}[section]
\newtheorem{Rmk}{Remark}[section]
\newtheorem{Cor}{Corollary}[section]
\newtheorem{Prop}{Proposition}[section]
\newtheorem{Thm}{Theorem}[section]
\newtheorem*{Thm*}{Theorem}
\newtheorem*{Prop*}{Proposition}
\newtheorem*{Cor*}{Corollary}
\newtheorem*{Lem*}{Lemma}
\theoremstyle{definition}
\newtheorem{Def}{Definition}[section]
\begin{document}
\title{Volume entropy and the Gromov boundary of flat surfaces}
\author{Klaus Dankwart}

\begin{abstract}
We consider the volume entropy of closed flat surfaces of genus $g\geq 2$ and area 1. 
We show that a sequence of flat surfaces diverges in the moduli space if and only if the volume entropy converges to infinity. Equivalently the Hausdorff dimension of the Gromov boundary of the isometric universal cover tends to infinity.\\
Moreover, we estimate the entropy of a locally isometric branched covering of a flat surface by the entropy of base surface and the geometry of the covering map. 
 \end{abstract}
\maketitle
\section{Introduction}
In the study of Teichm\"{u}ller space of closed Riemann surfaces,  half-translation structures play a central role. A half-translation structure on a surface defines a flat metric, see Section \ref{sectflatmetric} for details. Those flat metrics are  more special than general ones as  considered by \cite{Troyanov2007}, \cite{Troyanov1986}.\\   
The geometry of such flat metrics and their relationship to the hyperbolic metric in the same conformal class has been studied by  \cite{Duchin2010}, \cite{Rafi2007}. We  are interested  in which way the large-scale geometry of a flat surface $S$ is related to the large-scale geometry of a hyperbolic surface. To this end, we normalize $S$ to total area 1 and consider the locally isometric (or flat) universal cover $\pi:\tilde{S}\to S$. By the \v{S}varc-Milnor Lemma, $\tilde{S}$ is quasi-isometric to the Poincar\'{e} disc and so admits the following invariants, see Section \ref{sectgrhyp} for precise definitions:  
\begin{enumerate}[(1)]
 \item The volume entropy of $S$.
  \item  $\delta_{\inf}>0$, the infimum of all constants $\delta$ that  $\tilde{S}$ is $\delta$-hyperbolic in the sense of  Gromov. 
\item The Hausdorff dimension of the Gromov boundary of $\tilde{S}$.
 \end{enumerate}
On a hyperbolic surface of finite volume, these quantities are independent of the surface. This does not hold in the case of flat surfaces. Denote by $ \mathcal{Q}_{g}$ the moduli space of closed flat surfaces of genus $g\geq2$ and area 1 which can be identified with the unit cotangent space of the moduli space of Riemann surfaces.
  \begin{Thm*}[Theorem \ref{Thmentromodspace}, Corollary \ref{Corhdimmodspace} ]
For a sequence of flat surfaces  $S_i \in \mathcal{Q}_{g}$ the following conditions are equivalent.
\begin{enumerate}
 \item  $S_i$ diverges in $\mathcal{Q}_g$.
\item The volume entropy of $S_i$ tends to infinity. 
\item The Hausdorff dimension of the Gromov boundary of the locally isometric universal cover of $S_i$ tends to infinity. 
\end{enumerate}
\end{Thm*}
In the second part we estimate the volume entropy of a \textit{flat branched covering} $\pi: T \to S$,  i.e. a topological covering of flat surfaces that outside the branch points is a local isometry. Denote by $l_b(T)$ the minimal distance between branch points and by $\lambda(T)$ a measure of the combinatorics i.e.  the sum of the number of sheets and the number of branch points. 
 \begin{Thm*}[Theorem \ref{thmrcov}]
There is some constant $C>0$ and a function $a(S)>0$ that for each flat branched covering $\pi:T \to S$ the
 volume entropy $e(T)$ is bounded by the inequality 
 $$e(T)\leq (e(S)+1) \left(a(S)+\frac{C \cdot \log(\lambda(T))} {l_b(T)}\right)$$ 
The same holds for the Hausdorff dimension of the Gromov boundary.
 \end{Thm*}
 Moreover, we construct a family of examples which show that the bounds are asymptotically sharp.\\
 The paper is organized as follows. In Section \ref{sectprel} we recall the main known results about $\delta$-hyperbolic spaces and flat surfaces needed later on. In Section \ref{secttechlem} we establish some technical results on the geometry of flat surfaces and geodesics which  are known by the experts but did not appear in the literature.  In Section \ref{secthdimmodspace} we investigate how the invariants depend on the point in moduli space. In Section \ref{sectbrcover} we estimate the invariants for flat branched coverings.\\

 \textbf{Acknowledgement} I would like to thank my advisor Ursula Hamenst\"{a}dt for her support. I also like to thank my colleagues Sebastian Hensel and Emanuel Nipper for countless discussion. Parts of this work were done during my stay at the Universit\'{e} Paul Cezanne. I would a like to thank Pascal Hubert for inviting me and for many advices. This research was supported by the Bonn International Graduate School in Mathematics. 
 \section{Preliminaries about Gromov hyperbolic spaces and flat surfaces}\label{sectprel}
\subsection{Gromov hyperbolic spaces and their boundary}\label{sectgrhyp}
We recall the standard facts about proper $\delta$-hyperbolic spaces, compare \cite[ Chapter III] {BridsH1999}.\\ 
\textbf{Convention:} Any metric space is assumed to be complete, proper and geodesic.\\

A metric space $X$ is $\delta$-\textit{hyperbolic} if every geodesic triangle in $X$ with sides $a, b, c$ is $\delta$-\textit{slim}: The side $a$ is contained in the $\delta$-neighborhood of $b \cup c$.\\
$X$ admits a metric boundary which is defined as follows. Fix a point $p\in X$ and for two points $x,y\in X$ we define the \textit{Gromov product} $(x,y)_p :=\frac{1}{2}(d(x,p)+d(y,p)-d(x,y))$. We call a sequence $x_i$ \textit{admissible} if $(x_i,x_j)_p\rightarrow \infty$. We define two admissible sequences
 $x_i,y_i\subset X$ to be equivalent if $(x_i,y_i)_p \rightarrow \infty$. Since $X$ is hyperbolic, this defines an equivalence relation. The boundary $\partial X$ of $X$ is then the set of equivalence classes.\\ 
 The \textit{Gromov product on the boundary} is then 
$$(\eta,\zeta)_p=\sup\{\liminf\limits_{i,j}(x_i,y_j)_p~|~\{x_i\} \in \eta ,~\{y_j\} \in \zeta\} $$
\begin{Prop} \label{Prpgrmetric}
 Let $X$ be a $\delta$-hyperbolic space and let $\delta_{\inf}$ be the infimum of all Gromov hyperbolic constants. 
Moreover, let $\xi$ be defined by $2\delta_{\inf}\cdot \log(\xi)=\log(2)$. For any $p\in X$ there is a metric $d_{p,\infty}$ on $\partial X$ and a constant $\epsilon(\delta_{\inf}) < 1$ which satisfies:
$$\xi^{-(\eta, \zeta)_p}\geq d_{p,\infty}(\eta, \zeta)\geq (1-\epsilon(\delta_{\inf}))\xi^{-(\eta, \zeta)_p} $$ 
\end{Prop}
$d_{p,\infty}$ is a \textit{Gromov metric} and $(\partial X,d_{p,\infty}) $ the \textit{Gromov boundary}. A quasi-isometry between Gromov $\delta$-hyperbolic spaces extends to a homeomorphism between the Gromov boundaries.\\ 
Denote by $d_{\infty}$ the bilipschitz equivalence class of Gromov metrics $d_{p,\infty},~ p\in X$.

For a compact metric space $X$ and the locally isometric universal cover $\pi: \tilde{X} \to X$ let $\Gamma$ be the group of deck transformations. Fix $\tilde{p}\in \tilde{X}$ and let $N_{\tilde{p}}(R) :=|\tilde{q} \in \tilde{p}\Gamma | d(\tilde{p},\tilde{q})\leq R|$ be the \textit{counting function}. The \textit{volume entropy} is then 
$$e(X):= \limsup\limits_{R \to \infty}\frac{\log(N_{\tilde{p}}(R))}{R}$$
\textbf{Convention:} By entropy we mean volume entropy.\\

As $e(X)$ is independent of $\tilde{p}$, we skip the base point and abbreviate $N(R)$. \\
If the isometric universal cover $\tilde{X}$ is $\delta$-hyperbolic, \cite[Theorem 7.7]{Coorn1993} showed that for the bilipschitz equivalence class of metrics $d_{\infty}$ the entropy and the Hausdorff dimension of Gromov metric are related. 
\begin{Thm} \label{thmhdimentro}
If the universal cover $\tilde{X}$ is $\delta$-hyperbolic and $X$ is compact, let $d_{\infty}$ be the equivalence class of Gromov metrics on $\partial \tilde{X}$ induced by the constant $\xi$.\\ 
Then, the Hausdorff dimension of the Gromov boundary of $\tilde{X}$ is $\frac{e(X)}{\log(\xi)}$.
\end{Thm}
We remark that the Hausdorff dimension of the Gromov boundary of $\tilde{X}$ remains unchanged under scaling the metric on $X$. 
\subsection{Geometry of flat surfaces}\label{sectflatmetric}
We introduce the geometry of flat surfaces and refer to \cite{Minsk1992}, \cite{Streb1984}. \\
A \textit{half-translation structures} on a surface $X$ of genus $g\geq 2$ is a choice of charts that, away from a finite set of points $\Sigma$, the transition functions are half-translations. The pull-back of the flat metric in each chart gives a metric $d_\fl$ on $X-\Sigma$. We require that $d_\fl$ extends to a singular cone metric on $X$ with cone angle $k\pi$ in each point $\varsigma \in \Sigma$ where $k\geq 3$ is an integer. Then $S=(X,d_\fl)$ is a \textit{flat surface}. We scale the flat metric to area 1. \\
A \textit{straight line segment} on $S-\Sigma$ is defined as the pull-back of a straight line segment on $\R^2=\C$ in each chart. A straight line segment which emanates from one singularity and ends at another is a \textit{saddle connection}.\\
The flat metric admits an unoriented \textit{flat angle} as follows. A \textit{standard} neighborhood $U$ of a point $p \in S$ is isometric to a cone around $p$. On the boundary circle of $U$ we choose an orientation. Let $c_1,c_2$ be straight line segments, issuing from $x$. The complement $U -c_1 \cap c_2$ consists of two connected components $U_1,U_2$ which are isometric to euclidean circle sectors with angle $\vartheta_i,(i=1,2)$ possibly greater than $2\pi$. Choose $U_1$ that the arc on $\partial U$ which connects $c_1$ with $c_2$ in direction of the boundary orientation is on the boundary of $U_1$. Then $\angle_p(c_1,c_2)$ is the sector angle $\vartheta_1$. 
\begin{Lem} 
A path $c:[0,T]\to S$ is a local geodesic if and only if it is continuous and a sequence of straight lines segments outside $\Sigma$. In the singularities $\varsigma=c(t)$ the consecutive line segments make flat angle at least $\pi$ with respect to both boundary orientations. 
\end{Lem}
\begin{proof}
 We refer to \cite{Streb1984}. 
\end{proof}
One of the main concepts is the \textit{Gauss-Bonnet formula}, see \cite{Hubba2006}. Let $P$ be a compact flat surface with piecewise geodesic boundary and denote by $\chi(P)$ the Euler characteristic of $P$. For $x \in \interior{P}$, define $\vartheta(x)$ the cone angle of $x$ in $P$ and for $x\in \partial P$, $\vartheta(x)$ is the cone angle at $x$ inside $P$. Then 
$$2\pi\chi(P)=\sum_{x \in \interior{P}} (2\pi-\vartheta(x))+\sum \limits_{x\in \partial P} (\pi-\vartheta(x))$$ 
As a corollary, a flat surface does not contain geodesic bigons.
\begin{Prop}
 In any homotopy class of arcs with fixed endpoints on a closed flat surface there exists a unique local geodesic which is length-minimizing.\\ 
Moreover, in any free homotopy class of closed curves there is a length-minimizing locally geodesic representative.
\end{Prop}
\begin{proof}
In both cases, the existence follows from a standard Arzel\`{a}-Ascoli argument. The uniqueness follows from the absence of geodesic bigons. 
\end{proof}
A \textit{flat cylinder} of height $h$ and circumference $c$ on a flat surface $S$ is an isometric embedding of $[0,c] \times (0,h)/\sim, (0,t)\sim (c,t) $ into $S$. A flat cylinder is \textit{maximal} if it cannot be extended. \\ 
Like in hyperbolic geometry, closed local geodesics on flat surfaces do not have arbitrary intersections. 
\begin{Lem} \label{Lemintnumbergeod}
 Let $\alpha, \beta$ be closed curves on a flat surface $S$ and $\alpha_{\fl}, \beta_\fl$ be a choice of locally geodesic representatives in the free homotopy class. If the number of intersection points of $\alpha_\fl$, $\beta_\fl$ is bigger than the geometric intersection number $i([\alpha],[\beta])$, then the local geodesics $\alpha_\fl$ and $\beta_\fl$ share some arcs which begin and end at singularities. The arcs might be degenerated to singular points.
\end{Lem}
 \begin{proof}
If $\alpha_\fl$ and $\beta_\fl$ have more points in common than $i(\alpha, \beta)$, by absence of geodesic bigons, $\alpha_\fl$,$\beta_\fl$ share some arc. As local geodesics in flat surfaces are straight line segments outside the singularities, two local geodesics having some arc in common can ones drift apart at singularities. So they can ones share arcs with singularities as start- and endpoints. 
\end{proof}

\section{High-cylinder decomposition, short curve systems and geodesics joining saddle connections}\label{secttechlem}
We introduce some basic concepts concerning geodesics on a closed flat surface $S$ of genus $g \geq 2$ and area 1. Most are well-known but did not appear in the literature in this form. We summarize the results:
\begin{enumerate}[(1)]
 \item $S$ admits a decomposition $(\bigcup S_j,\bigcup C_i^*)$ into a disjoint union of subsurfaces $S_j$ of uniformly bounded diameter with geodesic boundary and high flat cylinders $C_i^*$ which connect the components $S_j$. Each free homotopy class of closed curves in a component $S_j$ contains a locally geodesic representative in $S_j$.
\item For any closed flat surface $S$ there is a constant $c>0$ which only depends on the topology of $S$, and there is a collection of essential closed curves $\alpha_{i,\fl}$ of length at most $c$ whose homotopy classes define a pants decomposition of $S$.
\item There is a constant $C_l(S)>0$ so that following holds:\\
For any two parametrized saddle connections $s,s'$ on $S$ there is a local geodesic $g$ which first passes through $s$ and eventually passes through $s'$ and which is of length at most $C_l(S)+l(c)+l(c')$. 
\end{enumerate}
\subsection{High-cylinder decomposition and short curve systems}\label{secthighcyl}
We recall the construction as in \cite[Section 5]{MasurS1991} to show that a flat surface of genus $g\geq 2$ and area $1$ consists of disjoint subsurfaces of uniformly bounded diameter and high flat cylinders which connect the components. The following two lemmas are well-known. 
\begin{Lem}\label{lemlargedistflatcyl}
 Let $\pi:\tilde{S} \to S$ be the flat universal cover of a flat surface. Let $\tilde{D} \subset \tilde{S}$ be a euclidean disc of radius $r$ i.e. an embedded disc not containing any singular point and let $\tilde{D}'$ be the subdisc of radius $r/2$ and the same center $\tilde{x}$.\\ 
If the restriction to $\tilde{D}'$ of the projection $\pi$ is not an embedding into $S$, then $\pi(\tilde{x})$ lies in a flat cylinder, which has circumference at most $r$ and height at least $\sqrt{\frac{3}{4}}r$. 
\end{Lem}
\begin{proof}
We sketch the proof. Since $\tilde{D}'$ does not contain a singularity and does not embed, there are two points in $\tilde{D}'$ which projects to the same point in $S$. Take two closest such points and let $\tilde{c}$ be the connecting geodesic in $\tilde{D}'$. $\tilde{c}$ does not pass through a singularity and therefore does not change direction. So it projects to a simple closed local geodesic. As $\tilde{D}$ is a euclidean disc, one can transport $\tilde{c}$ in direction of the normal bundle. One can proceed at least until one reaches the boundary of $\tilde{D}$. The projection of this family of curves is a flat cylinder of height at least $\sqrt{\frac{3}{4}}r$. 
\end{proof}
One deduces:
\begin{Lem}\label{lemlargedistflatcyl2}
There exists some constant $c_{\mathit{height}}>0$ such that for any closed flat surface $S$ of area $1$ and singularities $\Sigma$ the following holds: 
\begin{enumerate}
 \item Each point $x \in S$ with $d(x, \Sigma)>\sqrt{\frac{4}{3}} c_{\mathit{height}}$ is contained in a maximal flat cylinder of height at least $c_{\mathit{height}}$. 
\item Any two maximal flat cylinders of height at least $c_{\mathit{height}}$ are either disjoint or equal. 
\end{enumerate}
\end{Lem}
For a closed flat surface $S$ of area 1 and genus $g\geq 2$ let $\bigcup C_i$ be the disjoint union of all maximal flat cylinders $C_i$ each of height at least $c_{\mathit{height}}$. By Lemma \ref{lemlargedistflatcyl2}, $\bigcup C_i$ contains all points which are of distance at least $\sqrt{\frac{4}{3}} c_{\mathit{height}}$ to each singularity. Since the core curves of different maximal flat cylinders are simple, disjoint and pairwise not freely homotopic, the number of such cylinders is bounded from above by a constant which only depends on the topology of $S$.\\ 
For each such cylinder $C_i$ choose $C^*_i \subset C_i$ the closed central subcylinder of $C_i$ so that $C_i-C_i^*$ consists of two flat cylinders both of height $\frac{c_{\mathit{height}}}{3}$. Let $\bigcup S_j:= S-\bigcup_i C_i^* $ be the complement of the cylinders. 
\begin{Def} \label{Def}
We call the decomposition $\bigcup S_j$ of a closed flat surface $S$ of area 1 and genus $g \geq 2$ the high-cylinder decomposition. 
 \end{Def}
Each subsurface $S_j$ is of negative Euler characteristic and the boundary of $S_j$ is a disjoint union of geodesic core curves of some cylinders $C_i$ and especially disjoint from the singularities.\\ 
Since each point in $S_j$ is of uniformly bounded distance to a singularity and the number of singularities in $S$ is also uniformly bounded, the diameter of $S_j$ is bounded from above by some constant $c_{\mathit{diam}}$ which only depends on the topology of $S$. 
\begin{Lem}\label{Lemconvexsubsurf}
Let $S_j \subset S$ be a component of the high-cylinder decomposition of a closed flat surface $S$ and let $\alpha$ be some essential closed curve in $S_j$ which might be freely homotopic to a boundary component.\\ 
Then there is a locally geodesic representative $\alpha_\fl \subset S_j$ in the free homotopy class of $\alpha$ which is contained in $S_j$.
\end{Lem}
\begin{proof}
Assume first that $\alpha$ is freely homotopic to the multiple of some boundary component. Each boundary component of $S_j$ is the core curve of some cylinder $C_i$ and there exists a geodesic core curve of $C_i$ in $S_j$.\\ 
Assume next that $\alpha$ is not freely homotopic to the boundary and so cannot not be realized disjointly from $S_j$. It remains to show that the geodesic representative $\alpha_\fl$ of $\alpha$ does not intersect the boundary of $S_j$. But each boundary component of $S_j$ is the locally geodesic core curve of some cylinder $C_i $ and so, by Lemma \ref{Lemintnumbergeod} disjoint from $\alpha_\fl$.
 \end{proof}
\begin{Lem}\label{Lemshortcurvsubsurf}
Let $\bigcup S_j$ be the high-cylinder decomposition a flat surface $S$ of genus $g\geq 2$ and area 1 and let $c_{\mathit{height}}$ be the constant as in Lemma \ref{lemlargedistflatcyl2}.\\ 
Then, for each closed curve $\alpha \subset S$ of length less than $c_{\mathit{height}}$ there exists a length minimizing representative $\alpha_\fl$ of the free homotopy class $[\alpha]$ which is a contained in some subsurface $S_j$. 
\end{Lem}
\begin{proof}
 We can assume that $\alpha$ is length-minimizing in its free homotopy class of closed curves and therefore locally geodesic. Let $C_j$ be a maximal flat cylinder which contains a component of the complement $S-\bigcup S_j$ and let $\beta$ be a geodesic core curve of $C_j$. If $\alpha$ intersects $\beta$ then it has to cross through the flat cylinder $C_j$ and has length at least $c_{\mathit{height}}$, which is impossible. Therefore, $\alpha$ is freely homotopic to some curve which is disjoint from the core curves of all cylinders and can be homotoped in some subsurface $S_j$.\\ 
By Lemma \ref{Lemconvexsubsurf} there exists a geodesic representative $\alpha_{\fl} $ in the free homotopy class of $\alpha$ which is contained in $S_j$. 
\end{proof}

Moreover, we recall the existence of a maximal collection of essential closed curves on $S$ which are in the free homotopy class of disjoint simple closed curves.
\begin{Prop}\label{propshortpants}
For any closed flat surface $S=(X,d_\fl)$ there is a constant $c>0$, which only depends on the topology of $S$, and there is a collection of essential closed curves $\alpha_{i,\fl}$ of length at most $c$ whose homotopy classes define a pants decomposition of $S$.
 \end{Prop}
\begin{proof}
For the hyperbolic metric $ \sigma$ on $X$ in the same conformal class as $d_\fl$ it was shown by \cite{Bers1985} that there exists a constant $c_B$, which only depends on the topology of $X$, and such a maximal collection of essential disjoint simple closed curves $\alpha_{i,\sigma}$ of hyperbolic length at most $c_B$. \\
Let $c:=\frac{1}{2} c_B \exp\left(\frac{c_B}{2}\right)$ and let $\alpha_{i,\fl}$ be the flat geodesic representative in the free homotopy class of $\alpha_{i,\sigma}$. By \cite{Maskit1985}, the flat length of $\alpha_{i,\fl}$ is at most $c$.
 \end{proof}
\subsection{Geodesic connection of saddle connections}\label{sectconcasc}
In this part we show that for any two saddle connections on a closed flat surface $S$ one finds a connecting local geodesic of uniformly bounded length, so that the concatenation is a local geodesic. The discussion is motivated by \cite{Duchin2010}.
 \begin{Prop}\label{Propconsc}
 For a closed flat surface $S$ there is a constant $C_l(S)>0$ so that following holds:\\
For any two parametrized saddle connections $s,s'$ on $S$ there is a local geodesic $g$ which first passes through $s$ and eventually passes through $s'$ and which is of length at most $C_l(S)+l(c)+l(c')$. 
\end{Prop}
To show the proposition we need some additional lemmas.
For the flat universal cover $\pi:\tilde{S} \to S$ fix a point $\tilde{x}\in \tilde{S}$ and a set $\tilde{U}\subset S$. A point $\eta \in \partial \tilde{S}$ in the Gromov boundary is contained in the \textit {boundary shadow} $\partial sh_{\tilde{x}}(\tilde{U})\subset \partial \tilde{S}$ if and only if there is a representing sequence $\{\tilde{x}_i\}$ so that for almost all $i$ the geodesic connecting $\tilde{x}$ and $\tilde{x}_i$ intersects $\tilde{U}$. It $\eta$ is in the interior of $\partial sh_{\tilde{x}}(\tilde{U})$ this holds for all representatives. \\ 
It is well-known, see \cite[Chapter III]{BridsH1999} that the boundary shadow of an open set $\tilde{U}$ is an open set on the boundary. 
 \begin{Lem}
 Fix points $\tilde{x}\not= \tilde{\varsigma}$ on the flat universal cover of a flat surface so that $\tilde{\varsigma}$ is a singularity. Then $\partial sh_{\tilde{x}}(\{\tilde{\varsigma}\})$ contains an open subset of the Gromov boundary.
 \end{Lem}
 \begin{proof}
 Let $[\tilde{x}, \tilde{\varsigma}]$ be the geodesic connecting $\tilde{x}$ with $\tilde{\varsigma}$. As the cone angle at $\tilde{\varsigma}$ is at least $3\pi$, in a standard neighborhood of $\tilde{\varsigma}$ one finds a small open subset $\tilde{U}$ so that for any $\tilde{y} \in \tilde{U}$ the geodesic which connects $\tilde{x}$ and $\tilde{y}$ passes through $\tilde{\varsigma}$. Observe that $\partial sh_{\tilde{x}}(\tilde{\varsigma})\supset \partial sh_{\tilde{x}}(\tilde{U})$.
\end{proof}
A bi-infinite geodesic $\tilde{c}$ on a flat metric in the disc $\tilde{S}$ decomposes $\tilde{S}$ in two components $\tilde{S}^{\pm}$. 
The \textit{flat angle} $\angle_{\tilde{c}(t)}^+(\tilde{c})$ \textit{at the side} $\tilde{S}^+$ is the sector angle of the intersection of a standard neighborhood of $\tilde{c}(t)$ with $\tilde{S}^+$ 
\begin{Lem}\label{Lemcritcorecurve}
Let $\alpha$ be a closed local geodesic on a closed flat surface $S$ and let $\tilde{\alpha}$ be a complete lift of $\alpha$ in the flat universal cover. $\alpha$ is freely homotopic to a core curve of a flat cylinder if and only if at one side $\tilde{S}^{+}$ the flat angle $\angle^+_{\tilde{\alpha}(t)}(\tilde{\alpha})$ is $\pi$ for all $t$.
\end{Lem} 
\begin{proof}
We refer to \cite[Lemma 20]{Duchin2010}.
\end{proof}
\begin{Lem}\label{lemconnscclosed}
 Let $s$ be a parametrized saddle connection on a closed flat surface $S$ and let $\alpha$ be a parametrized closed local geodesic on $S$ which is not freely homotopic to a core of flat cylinder. \\
Then there exists a local geodesic $g$ which first passes through $s$ and eventually through $\alpha$.
\end{Lem}
\begin{proof}
Let $\tilde{s}$ be a parametrized lift of $s$ to the flat universal cover $\tilde{S}$ and let $\tilde{p}$ resp. $\tilde{\varsigma}$ be the starting point resp. endpoint of $\tilde{s}$. \\
As the visual metric boundary of the Poincar\'{e} disc and the Gromov boundary of $\tilde{S}$ are homeomorphic and the Deck transformation group $\Gamma$ acts on both boundaries in same topological way there is a parametrized complete lift $\tilde{\alpha}$ of $\alpha$ which is disjoint from $\tilde{s}$ and whose positive endpoint is contained in the interior of the boundary shadow $\partial sh_{\tilde{p}}(\{\tilde{\varsigma}\})$. \\
Let $\tilde{S}^+$ be the component of $\tilde{S}-\tilde{\alpha}$ that contains $\tilde{s}$ and let $\tilde{c}_t:[0,a_t] \to \tilde{S}$ be the geodesic which connects $\tilde{p}$ with $\tilde{\alpha}(t)$. There is a time $t_0$ that for all $t\geq t_0$, $ \tilde{c}_t$ passes through $\tilde{\varsigma}$ and that $\tilde{c}_{t_0}$ only shares its endpoint with $\tilde{\alpha}$.\\
As $\tilde{\alpha}$ is not in the free homotopy class of a core curve of a flat cylinder, by Lemma \ref{Lemcritcorecurve} there is a sequence of times $a\in \{a_0+n\cdot l_\fl(\alpha),~n\in \Z \} $ where the flat angle $\angle_{\tilde{\alpha}(a)}^+(\tilde{\alpha})$ at the side in $\tilde{S}^+$ is fixed $\pi+\epsilon$ for some $\epsilon>0$.\\ 
 So, there is some $t_1>t_0$ that the subarc $\left. \alpha \right|_{[t_0,t_1]}$ has excessive flat angle at least $\pi$:
$$\sum\limits_{t_0<t<t_1} \left(\angle_{\tilde{\alpha}(t)}^+(\tilde{\alpha})-\pi\right)> \pi$$
By the Gauss-Bonnet formulae, no geodesic triangle in $\tilde{S}^+$ whose three sides only intersect at the vertices contains the whole side $\left. \alpha \right|_{[t_0,t_1]}$.\\ 
Consider the triangle with sides $\tilde{c}_{t_0},\tilde{c}_{t_1},\left. \alpha \right|_{[t_0,t_1]}$. Since $\left. \alpha \right|_{[t_0,t_1]}$, $\tilde{c}_{t_0}$ only intersect at their endpoint, $\left. \alpha \right|_{[t_0,t_1]}$ has to coincide with $\tilde{c}_{t_1}$ for some time.\\ 
 The projection $g:=\pi(\tilde{c}_{t_1})$ is a local geodesic which first passes through $s$ and eventually through a piece of $\alpha$ and so can be extended to a local geodesic which passes through the whole curve $\alpha$.
\end{proof}
\begin{Cor}
For a closed flat surface $S$ and for any two parametrized saddle connections $s_1,s_2$ on $S$ there is a local geodesic $g$ which first passes through $s$ and eventually passes through $s'$.
 \end{Cor}
\begin{proof}
 Fix a parametrized closed local geodesic $\alpha$ which is not freely homotopic to a core curve of a flat cylinder. By Lemma \ref{lemconnscclosed} there are local geodesics $g_1$ resp. $g_2$ which first pass through $s_1$ resp. $s_2^{-1}$ and eventually through $\alpha$ resp. $\alpha^{-1}$. The concatenation $g=g_1 * g_2^{-1}$ is a local geodesic with the required properties. 
\end{proof}
It remains to show that we can choose $g$ of uniformly bounded length.
\begin{Lem}\label{lemscdense}
 Fix a point $x \in S$ on a closed flat surface $S$ of genus $g\geq 2$ and $\theta$ an outgoing direction at $x$. Let $\angle_\varsigma$ be the flat angle at $\varsigma$ with respect to a choice of orientation.\\ 
Then for any $\epsilon$ there exists a geodesic $c$ which connects $x$ with a singularity so that $\angle_x(c, \theta)\leq \epsilon$. 
\end{Lem}
\begin{proof}
This is well-known. We refer to \cite[Proposition 3.1]{Vorob1996}.
 \end{proof}
\begin{Lem} \label{lemsccover}
Let $S$ be a flat surface and let $\varsigma\in S$ be a singularity. There exist $4$ saddle connections $s_1, \ldots, s_4$ emanating from $\varsigma$ with the following property:\\ 
Let $c$ be a local geodesic with endpoint $\varsigma$. The concatenation of $c$ with at least one $s_i$ is again a local geodesic.
 \end{Lem}
\begin{proof}
This follows from the fact that the set of directions of saddle connections is dense and that the cone angle at each singularity is at least $3\pi$ 
\end{proof}
\begin{proof}[Proof of the Proposition]
Choose a set of parametrized saddle connections $s_{i}$ on the flat surface $S$ with property that any saddle connection $s$ can be concatenated with some saddle connection $s_{i}$ to an extended local geodesic.\\ 
By Lemma \ref{lemsccover} it suffices to choose a finite number of such saddle connections. For each pair $s_{i},s_{j}$, there is a local geodesic $g_{i,j}$ which first passes through $s_i$ and eventually through $s_j$. Since there are only finitely many pairs, the length of $g_{i,j}$ is bounded from above by a constant $ C_l(S)$. For any two saddle connections $s,s'$ there are two saddle connections $s_i$, $s_j$ such that the concatenation $s* s_i$ and $s_j* s'$ is a local geodesic. So $g:=s*g_{ij}* s'$ is also a local geodesic of length 
 $$l_\fl(s*g_{ij} *s)\leq C_l(S)+l(s)+l(s')$$ 
 \end{proof}
\section{Hausdorff dimension and entropy in moduli space}\label{secthdimmodspace}
Let $S$ be a closed flat surface and $\pi:\tilde{S}\to S$ the flat universal cover which is a $\delta$-hyperbolic geodesic metric space. We defined $\delta_{\inf}(\tilde{S})$ the infimum of all $\delta'$ so that $\tilde{S}$ is $\delta'$-hyperbolic.\\
To estimate $\delta_{\inf}$ explicitly in geometric terms of $\tilde{S}$, let $\Sigma_S$ resp. $\Sigma_{\tilde{S}}$ be the singularities on $S$ resp. $\tilde{S}$ and denote by $\rho(S)$ the \textit{packing density}
$$\rho(S):=\sup_{x \in S} d_\fl(x,\Sigma_S)=\sup_{\tilde{x} \in \tilde{S}} d_\fl(\tilde{x},\Sigma_{\tilde{S}})<\infty$$
\begin{Prop} \label{PSgromov}
 The flat universal cover $\pi: \tilde{S} \to S$ of a closed flat surface is $\delta$-hyperbolic with $\rho(S)/2 \leq \delta_{\inf} \leq 2\rho(S)$. If $S$ has area $1$, there is a lower bound on $\delta_{\inf}(S)$ which only depends on the topology of $S$. 
\end{Prop}
\begin{proof}
Let $\triangle(\tilde{x}_1, \tilde{x}_2, \tilde{x}_3) $ be a triangle in $\tilde{S}$ with vertices $\tilde{x}_i$. 
\begin{figure}[htbp]
 \centering
 \fbox{
 \includegraphics[width=0.4\textwidth]{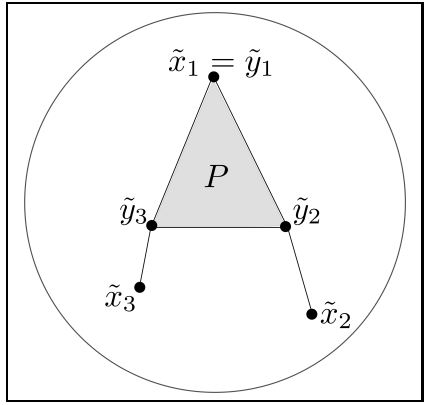}
 }
 \caption{The geodesics which correspond to the triangle of the $\tilde{x}_i$ might share some arc until they spread apart. We call the spread point $\tilde{y}_i$.}
 \label{figdeltahyp}
\end{figure}
 The sides emanating from a point $\tilde{x}_i$ might coincide for some time but after spreading, they remain disjoint. Denote by $\tilde{y}_i$ the point they start spreading apart. It suffices to show that the smaller triangle $\triangle(\tilde{y}_1, \tilde{y}_2, \tilde{y}_3)$ is $2\rho(S)$-slim.\\ 
Denote by $\tilde{a},\tilde{b},\tilde{c}$ the sides of $\triangle(\tilde{y}_1, \tilde{y}_2, \tilde{y}_3)$ and let $P$ be the interior of $\triangle(\tilde{y}_1, \tilde{y}_2, \tilde{y}_3)$ which is either empty or a topological disc, see Figure \ref{figdeltahyp}.\\ 
If $P=\emptyset$, the triangle is a tripod, consequently it is 0-slim, so assume that $P$ is not empty. 
By the Gauss-Bonnet formula, $P$ does not contain a singularity and so no disc of radius $\rho(S)$ is contained in $P$.\\
For a point $\tilde{p}$ on the side $\tilde{a}$ we have to show that the distance of $\tilde{p}$ to $\tilde{b} \cup \tilde{c}$ is at most $2\rho(S)$. As singular points are discrete, we can assume that $\tilde{p}$ is a regular point and not an endpoint of $\tilde{a}$. Choose a geodesic line segment $\tilde{c}$ that issues from $\tilde{p}$ perpendicular to $\tilde{a}$ inside $P$ and let $\tilde{x} \in \tilde{c}$ be a point with distance $\rho(S)$ to $\tilde{p}$. It follows from comparison geometry that the distance between $\tilde{x}$ and $\tilde{a}$ is also at least $\rho(S)$, so the open ball $\tilde{B}$ of radius $\rho(S)$ about $\tilde{x}$ intersects $P$ but does not intersect $\tilde{a}$. As $\tilde{B}$ is not a subset $P$, $\tilde{p}$ is of distance at most $2\rho(S)$ to $\tilde{b}\cup\tilde{c}$.\\

We have to show that there exists a triangle which is not $\rho(S)/2$ slim. Recall that by definition of the packing density there is a ball of radius $\rho(S)$ in $\tilde{S}$ which does not contain a singularity hence is isometric to a euclidean disc. Inscribe a maximal equilateral triangle which is not $\rho(S)/2$ slim.\\ 

Finally, we show that for a flat surface $S$ of area 1 there is a lower bound on $\rho(S)$ which only depends on the genus of $S$.\\ 
Let $\varsigma_i, 1\ldots n$ be the set of singularities in $S$ with cone angle $ n_i\pi $. By Gauss Bonnet formula 
$\sum (n_i-2)=2g-2$.\\
For $\epsilon>0$ around each singularity $\varsigma_i$ choose a disc in $S$ of radius 
$$r=\sqrt{\frac{1}{\sum_i n_i\pi}}-\epsilon$$ 
The area of the union of the discs is at most $\sum_i r^2 n_i\pi <1$. Therefore, there exists a point in the complement of the discs which has distance at least $r$ to each singularity and so $\rho(S)\geq r$. Since $n_i\geq 3$ for each $i$, the constant $r$ has a lower bound which only depends on the genus of $S$.
\end{proof}
We investigate the Hausdorff dimension of the Gromov boundary with respect to the bilipschitz equivalence class of Gromov metrics $d_{\infty}$. 
\begin{Lem} \label{Lemhdmlowerbound}
 Let $S$ be a closed flat surface and $\pi:\tilde{S}\to S$ the flat universal cover. Let $d_{\infty}$ be a Gromov metric on the boundary. The Hausdorff dimension of the boundary $(\partial\tilde{S},d_\infty)$ is at least $1$.
\end{Lem}
\begin{proof}
The Gromov boundary $\partial\tilde{S}$ is homeomorphic to the boundary of the Poincar\'{e} disc, hence it is a topological circle. This proves the statement as topological dimension is a lower bound for the Hausdorff dimension.
 \end{proof}
We investigate how the Hausdorff dimension and the entropy vary under slight changes of the flat surfaces. The \textit{moduli space of flat surfaces} $\mathcal{Q}_g$is the set of all isometry equivalence classes of flat surfaces of genus $g\geq 2$ and area $1$. For $S \in \mathcal{Q}_g$ and $\epsilon>0$ let 
$$B_S(\epsilon):=\{ S' \in \mathcal{Q}_g: \exists f:S' \to S,f \mbox{ is a }(1+\epsilon)-\mbox{bilipschitz homeomorphism}\}$$
The sets $B_S(\epsilon), S\in \mathcal{Q}_g, ~\epsilon>0$ form a basis of the topology on $\mathcal{Q}_g$.\\ 
We consider the following mappings 
\begin{enumerate}
\item $l_0: \mathcal{Q}_g \to \R$ the length of the shortest essential simple closed curve on the flat surface $S$.
\item $e: \mathcal{Q}_g \to \R$ the entropy of $S$. 
\item $\delta_{\inf}: \mathcal{Q}_g \to \R$, the minimal Gromov hyperbolic constant of the flat universal cover $\tilde{S}$. 
 \item $\mathit{hdim}: \mathcal{Q}_g \to \R$ the Hausdorff dimension of the Gromov boundary of the universal cover with respect to the Gromov metric $d_\infty$.
 \end{enumerate}
\begin{Lem} \label{Propcontimodspace}
 The functions $l_0,e,\delta_{\inf}$ and $\mathit{hdim}$ are continuous in moduli space.
\end{Lem}
\begin{proof}
For $e$ and $l_0$ this is obvious. To observe, that $\delta_{inf}$ is continuous we refer to \cite[Chapter III]{BridsH1999}. Then, by Theorem \ref{thmhdimentro} $\mathit{hdim}$ is continuous as well.
\end{proof}
Consequently, in compact sets of moduli space the Hausdorff dimension of the Gromov boundary with respect to the metric $d_\infty$ is bounded from above. It remains to investigate the quantities under the degeneration of the flat surface.\\ 
A sequence $S_i \in \mathcal{Q}_g$ is \textit{divergent} if it eventually leaves every compact subset of $\mathcal{Q}_g$.\\
As $l_0: \mathcal{Q}_g \to \R_+$ is continuous, in compact subsets of $\mathcal{Q}_g$, $l_0$ has a positive lower bound. It is well-known, that the converse also holds. 
\begin{Prop}
For $\epsilon>0$ the set $l_0^{-1}((\epsilon, \infty)) \subset \mathcal{Q}_g$ is contained in a compact subset of $\mathcal{Q}_g$.
\end{Prop}
The proof mainly uses the Mumford compactness Theorem and results from extremal length.\\ 
To show that the entropy $e(S_i)$ tends to infinity if and only if the sequence $S_i\in \mathcal{Q}_g $ diverges, we need the following lemmas:
\begin{Lem}\label{Leminterscurves}
Let $S$ be a closed flat surface of genus $g\geq 2$ and let $\alpha, \beta:[0,1]\to S$ be closed curves sharing at least one point $p:=\alpha(0)=\beta(0)=\alpha(1)=\beta(1)$ so that the length $a=l_\fl(\alpha)$ is at most the length $b=l_\fl(\beta)$. \\ 
If the group $<\alpha, \beta>$, considered as a subgroup of $\pi_1(S,p)$, is neither cyclic nor trivial, then the entropy of $S$ can be estimated in terms of $a,b$: 
$$e(S)\geq \max\left\{\frac{\log(b)-\log(a)}{2b}, \frac{\log(2)}{b}\right\}$$ 
\end{Lem}
We emphasize that $\alpha, \beta$ are not necessarily local geodesics. Length in this context actually means the length of $\alpha, \beta$ and not the length of geodesic representatives. Moreover, we do not require that $\alpha$ and $\beta$ intersect transversely.
 \begin{proof}
The group $<\alpha, \beta><\pi_1(S,p)$ is free and as $<\alpha, \beta>$ is neither cyclic nor trivial, the positive semi-group $<\alpha, \beta>_{s} $ is free with respect to the generating system $\{\alpha, \beta\}$.\\ 
Let $\pi:\tilde{S}\to S$ be the flat universal cover and let $\tilde{p}\in \pi^{-1}(p)$ be a preimage of $p$. Choose connected arcs $\tilde{\alpha}, \tilde{\beta}:[0,1] \to \tilde{S}$, emanating from $\tilde{p}$, so that $\pi \circ \tilde{\alpha}=\alpha, \pi \circ \tilde{\beta}=\beta$. Let $\tilde{p}_\alpha:=\tilde{\alpha}(1)$, resp. $\tilde{p}_\beta:=\tilde{\beta}(1)$ be the endpoints. By definition of $a,b$, 
$$d_\fl(\tilde{p}_\alpha, \tilde{p})\leq a,~d_\fl(\tilde{p}_\beta, \tilde{p})\leq b$$
Let $\gamma_\alpha$ resp. $\gamma_\beta$ be the element of the Deck transformation group which maps $\tilde{p}$ to $\tilde{p}_\alpha$ resp. $\tilde{p}_\beta$ and let $\Phi$ be the canonical isomorphism of the positive semi-group of words with letters $\alpha, \beta$ to the semi-group $<\gamma_\alpha, \gamma_\beta>_{s}$ with distinguished generating system $\{\gamma_\alpha, \gamma_\beta\}$ which is defined as $\Phi(\alpha):=\gamma_\alpha, \Phi(\beta):=\gamma_\beta$.\\ 
 Let $w=a_1\ldots a_{k+l}$ be a word with letters in $\{\alpha,\beta \} $ which contains $k$ times the letter $\alpha$ and $l$ times the letter $\beta$. Let $w_i$ be the sub-word of $w$ truncated after the $i$-th letter.\\ 
Due to the triangle inequality
$$d_\fl(\tilde{p}, \Phi(w)(\tilde{p}))\leq \sum_{i=1}^{k+l} d_\fl(\Phi(w_{i-1})(\tilde{p}), \Phi(w_{i})(\tilde{p}))$$
Since $\Gamma$ is a group of isometries, it follows that the distance between $\tilde{p}$ and its image under $\Phi(w)$ can be estimated by the following formula: 
$$d(\tilde{p}, \Phi(w)(\tilde{p}))\leq ka+lb$$
Let $U(R )$ be the set of words so that the number of $\alpha$-letters is $\lfloor R/a\rfloor$ and the number of $\beta$-letters is $\lfloor R/b\rfloor$. The function $\lfloor * \rfloor: \R \to \N$ rounds down each number. 
 The cardinality of $U(R)$ is
\begin{eqnarray*}
|U(R)|&=& \binom{\lfloor R/a\rfloor+\lfloor R/b\rfloor} {\lfloor R/b\rfloor}\geq\frac{\prod\limits_{i=1}^{\lfloor R/b \rfloor} (i+R/a-1)}{\lfloor R/b \rfloor!} \geq \left(\frac{b}{a}-\frac{b}{R}\right)^{R/b-1}
\end{eqnarray*}
 For a word $w\in U(R)$ the distance $d(\tilde{p}, \Phi(w)\tilde{p})$ is at most $2R$. As the Deck transformation group acts freely we estimate the counting function
$$N_{\tilde{p}}(2R)\geq | U(R)| \geq \left(\frac{b}{a}-\frac{b}{R}\right)^{R/b-1}$$ 
Consequently,
$$e(S)\geq \frac{\log(b)-\log(a)}{2b}$$
The other inequality is analogous. Since $a\leq b$, for any word $w$ of length $n$ the distance can be estimated by 
$$d(\tilde{p}, \Phi(w)(\tilde{p}))\leq bn$$
Let $V(R)$ be the set of words of length $\lfloor R/b\rfloor$. The cardinality of $V(R)$ is $2^{\lfloor R/b\rfloor}$ and therefore
$$N_{\tilde{p}}(R)\geq |V(R)|\geq 2^{\lfloor R/b\rfloor}\geq \frac{1}{2} 2^{R/b}$$
 \end{proof}
 The following corollary is an immediate consequence. 
\begin{Cor}\label{Corshortcurvlargeentro}
Let $S_i$ be a sequence of flat surfaces and let $\alpha_i, \beta_i$ be closed curves on $S_i$ satisfying the conditions as in Lemma \ref{Leminterscurves}. If the length of $\alpha_i$ tends to zero and the length of $\beta_i$ is bounded from above, then the entropy of $S_i$ tends to infinity. 
 \end{Cor}
For a point $S \in \mathcal{Q}_g$ which contains a short simple closed curve $\alpha$ it is the goal to show that there exists a curve $\beta$ whose length is uniformly bounded from above so that $\alpha, \beta$ satisfy the conditions of Lemma \ref{Leminterscurves}. We recall the high-cylinder decomposition:\\ 
There exists a finite disjoint union of open subsurfaces $\bigcup S_j \subset S$ of negative Euler characteristic so that the diameter of each subsurface $S_j$ is bounded from above by some constant $c_{\mathit{diam}}$ which only depends on the topology of $S$. The complement $S-\bigcup S_j$ is contained in a disjoint union of maximal flat cylinders $\bigcup C_i$ whose heights are bounded from below by a uniform constant $c_{\mathit{height}}$. If $C_i$ intersect some subsurface $S_j$ then $S_j$ contains a geodesic core curve of $C_i$.\\
For each closed curve $\alpha$ which is contained in a subsurface $S_j$, there exists a geodesic representative $\alpha_\fl$ in the free homotopy class of $\alpha$ which is also contained in $S_j$. Moreover, let $\alpha$ be a closed curve in $S$ whose length is less than $c_{\mathit{height}}$. Then there exists a geodesic representative $\alpha_\fl$ in the free homotopy class of $\alpha$, which is contained in some subsurface $S_j$.
 \begin{Lem}\label{lemcloseshortcurves}
For $S\in \mathcal{Q}_g$ let $\bigcup S_j$ be the high-cylinder decomposition of $S$ and let $\alpha$ be a simple closed curve in $S_1$ which might be freely homotopic to a boundary curve.\\
Then there exists some uniform constant $c>0$ which only depends on the genus of $S$ and a simple closed local geodesic $\beta_\fl \subset S$ which intersects $S_1$ so that the length of $\beta_\fl$ is bounded from above by $c$ and no multiple of $\beta_\fl$ is in the free homotopy class of any multiple of $\alpha$ in $S$.
\end{Lem}
\begin{proof}
By Proposition \ref{propshortpants} there is a maximal system of closed curves $\beta_{\fl,i}$ which can be homotoped to disjoint simple closed curves and the length of $\beta_{\fl,i}$ is bounded from above by a constant $c$ which only depends on the topology of $X$.\\
Assume first that there is a curve $\beta_{\fl,i}$ which intersects $S_1$ essentially, but cannot be homotoped inside of $S_1$. No multiple of $\beta_{\fl,i}$ is in the free homotopy class of a multiple of $\alpha$ and $\beta:=\beta_{\fl,i}$ also intersects $S_1$.\\
If such a curve does not exist, recall that the Euler characteristic of $S_1$ is negative and so there are two curves $\beta_{\fl,i_1}, \beta_{\fl,i_2}$ which can be homotoped inside of $S_1$, possibly as boundary components.\\
By Lemma \ref{Lemconvexsubsurf} we can assume that $\beta_{\fl,i_j},j=1,2$ are contained in $S_1$. For at least one $j$ no multiple of $\beta:=\beta_{\fl,i_j}$ is freely homotopic to a multiple of $\alpha$.
 \end{proof}
\begin{Thm} \label{Thmentromodspace}
A sequence $S_i \in \mathcal{Q}_g$ diverges if and only if the entropy $e(S_i)$ tends to infinity. 
\end{Thm}
\begin{proof}
As the entropy depends continuously on the flat surface, it remains to show that for a divergent sequence $S_i \in \mathcal{Q}_g$ the entropy tends to infinity. \\
There exists a sequence of essential simple closed curves $\alpha_i$ in $S_i$ so that the length of $\alpha_i$ tends to zero. Let 
$$\bigcup\limits_j S_{i,j}\subset S_i$$ 
be the high-cylinder decomposition of $S_i$. If $\alpha_i$ is shorter than $c_{\mathit{height}}$ the minimal height of the removed cylinders, by Lemma \ref{Lemshortcurvsubsurf}, $\alpha_i$ is contained in some subsurface $S_{i,j}$. By Lemma \ref{lemcloseshortcurves}, there exists an essential curve $\beta'_i$ of length less than some uniform constant $c$ which intersects $S_{i,j}$ so that no multiple of $\beta'_i$ is freely homotopic to a multiple of $\alpha_i$.\\ 
As the diameter of $S_{i,j}$ is bounded from above by some uniform constant $c_{\mathit{diam}}$, choose a short excursion from $\beta'_i$ to $\alpha_i$ of length at most $c_{\mathit{diam}}$, concatenate $\beta_i'$ with the excursion and call the resulting closed curve $\beta_i$.\\ 
The length of $\beta_i$ is uniformly bounded and the two curves $\alpha_i, \beta_i$ satisfy the conditions of Lemma \ref{Leminterscurves}, so by Corollary \ref{Corshortcurvlargeentro} the entropy of $e(S)$ tends to infinity.
\end{proof}
As a consequence we obtain
\begin{Cor}\label{Corhdimmodspace}
Let $S_i\in \mathcal{Q}_g$ be a sequence of flat surfaces. The Hausdorff dimension of the boundary $\mathit{hdim}(S_i)$, with respect to the Gromov metric $d_\infty$, tends to infinity if and only if the sequence $S_i$ diverges. 
\end{Cor}
\begin{proof}
By Theorem \ref{Thmentromodspace} $S_i$ diverges if and only if the entropy $e(S_i)$ tends to infinity. By Theorem \ref{thmhdimentro} 
 the Hausdorff dimension equals $\frac{e(S_i)}{\log(\xi(\delta_{\inf}(S_i)))}$ where $\log(\xi(\delta_{\inf}(S_i)))$ is bounded from above by some positive constant which only depends on the topology of $S$.
\end{proof}
 \section{Hausdorff dimension under branched coverings}\label{sectbrcover}
 Let $\pi:T \to S$ be a finite-sheeted flat branched covering of flat surfaces. The preimages of the singularities $\Sigma_S$ on $S$ and the set of branch points form the singularities $\Sigma_T$ on $T$ and the area of $T$ is the product of the number of sheets and the area of $S$.\\ 
We investigate how the entropy and the minimal Gromov hyperbolic constant of the flat universal covers of $S$ and $T$ are related.
\begin{Prop}
For an $n$-sheeted branched flat covering $\pi:T\to S$ with $k$ branch points let $\tilde{S}$ resp. $\tilde{T}$ be the flat universal cover of $S$ resp. $T$.\\ 
Then the minimal Gromov hyperbolic constant $\delta_{\inf}(\tilde{T})$ can be estimated in terms of $\delta_{\inf}(\tilde{S})$ and $k$. 
$$\frac{\delta_{\inf}(\tilde{S})}{24(k+1)}\leq \delta_{\inf}(T)\leq 4\delta_{\inf}(\tilde{S})$$ 
\end{Prop}
\begin{proof}
Recall that by Proposition \ref{PSgromov} the packing density $\rho(S)$ and $\delta_{\inf}(\tilde{S})$ are related by $\frac{\rho(S)}{2}\leq \delta_{\inf}(\tilde{S})\leq 2\rho(S)$.\\ 
We first show that the packing density of $T$ is at most the packing density of $S$. For a point $y \in T$ let $x=\pi(y) \in S$ be the image of $y$ and let $c$ be a shortest geodesic connecting $x$ with a singularity. Choose a lift $c'$ of $c$ to $T$ which emanates from $y$ and which is maximal with respect to the property that it does not contain a branch point. Since branch points are singularities and since the flat covering is locally isometric, $d(y, \Sigma_T)\leq l(c')\leq l(c)=d(x,\Sigma_S)$. Consequently, 
$$\rho(T )\leq \rho(S)$$ 
On the other hand, let $B_T\subset T$ be the set of all branch points on $T$ and let $B_S:= \pi(B_T) \subset S$ be the image of $B_T$ which contains at most $k$ points. Let $B_{\tilde{S}}:=\pi^{-1}(B_S) \subset \tilde{S}$ be the preimage of $B_S$ in the flat universal cover and let $\tilde{D}_1$ be a disc with center $\tilde{x}$ of radius $\rho(S)$ in $\tilde{S}$ which does not contain a singularity. 
 Let $\tilde{D}_2 \subset\tilde{D}_1$ be the euclidean sub-disc of radius $\frac{\rho(S)}{2}$ and center $\tilde{x}$. We distinguish two cases.
 \begin{enumerate}
 \item The projection of $\tilde{D}_2$ embeds into $S$.\\ 
Then there are at most $k$ points of $B_{\tilde{S}}$ in $\tilde{D}_2$ and so the discs of radius $r=\frac{\rho(S)}{3\sqrt{k}}$ around each point of $B_{\tilde{S}}\cap \tilde{D}_2$ do not cover $\tilde{D}_3\subset \tilde{D}_2$, the subdisc of center $\tilde{x}$ and radius $\frac{\rho(S)}{3}$.\\ 
 Let $\tilde{x}' \in \tilde{D}_3$ be a point which is of distance at least $\frac{\rho(S)}{3\sqrt{k}}$ to any point in $B_{\tilde{S}}\cap \tilde{D}_2$ and of distance $\frac{1}{6}\rho(S)$ to the boundary of $\tilde{D}_2$. So, it has distance at least $\frac{\rho(S)}{6\sqrt{k}}$ to any point in $B_{\tilde{S}}$ and any singularity.\\ 
Let $y'\in T$ be a point which projects to the same point in $S$ as $\tilde{x}'$. $y'$ is also of distance at least $\frac{\rho(S)}{6\sqrt{k}}$ to the singularities $\Sigma_T \subset T$. 
 Consequently,
 $$\rho(T) \geq d (y', \Sigma_T)\geq \frac{\rho(S)}{6\sqrt{k}}\geq \frac{\rho(S)}{6(k+1)}$$ 
 \item $\tilde{D}_2$ is immersed but not embedded. By Lemma \ref{lemlargedistflatcyl}, the projection of $\tilde{D}_2$ in $S$ is contained in a cylinder $C'$ of height at least $\frac{\sqrt{3}\rho(S)}{2}$. As the cardinality of $B_S$ is at most $k$, there is a sub-cylinder $C\subset C'$ of height $\frac{\sqrt{3}\rho(S)}{2(k+1)}$ which does not contain an element of $B_S$. Therefore, one finds a point $x' \in C$ of distance $\frac{\sqrt{3}\rho(S)}{4(k+1)}$ to a singularity and to $B_S$. 
 \end{enumerate}
In both cases 
$$\rho(T)\geq \frac{\rho(S)}{6(k+1)}$$
Using the relationship between $\rho$ and $\delta_{\inf}$ one estimates: 
 $$\frac{\delta_{\inf}(\tilde{S})}{24(k+1)}\leq \delta_{\inf}(\tilde{T}) \leq 4\delta_{\inf}(\tilde{S})$$ 
\end{proof}
We construct a family of examples which show that the bounds are asymptotically sharp:\\ 
A \textit{Strebel differential with one cylinder} is a flat surface which, after removing a finite union of saddle connections, is isometric to a single flat cylinder. Strebel differentials with one cylinder is dense in the moduli space of flat structures \cite{Masur1979}.\\ 
Let $S'$ be a Strebel differential with one cylinder of area $1$ and assume that the core curve of the cylinder is horizontal. One can stretch the vertical component by a large factor $\lambda$ and shrink the horizontal component by $\lambda^{-1}$. The resulting flat surface $S$ is again a Strebel differential with one cylinder of area $1$. The horizontal cylinder is of short circumference and large height $h$. The Gromov hyperbolic constant of the flat universal cover $\delta_{\inf}(\tilde{S})$ is nearly $\frac{h}{2}$. We can distribute $k$ points $x_i$ on $S$ so that each point in $S$ has distance at most $\frac{2h}{k+1}$ to some point $x_i$. We can construct a flat branched covering $\pi: T \to S$, so that each preimage of each point $x_i$ is a branch point and therefore a singularity. \\
In this special example one observes that the quantity $\delta_{\inf}(\tilde{T})$ is nearly $\frac{\delta_{\inf}(\tilde{S})}{k+1}$. The precise statement is summarized in the following remark: 
\begin{Rmk}
For any $g\geq 2$ there exists a family of branched covering $\pi_i:T_i \to S_i, S_i \in \mathcal{Q}_g, i\in \N$ so that the number of branch points on $T_i$ equals $i$ and $i\cdot \delta_{inf}(\tilde{T}_i)$ is bounded from above independent of $i$. 
\end{Rmk}
To compute the entropy of a flat covering we first consider the following example:\\ 
Let $S$ be a flat surface and $D \subset S$ be a small euclidean disc. On $D$ one can make a small slit of length $l$ and take $n\geq 3$ copies $S_1, \ldots S_n$ of $S$ endowed with the same slit. The copies can be isometrically glued along the slit. The right-hand side of the slit in $S_{i}$ is glued to the left-hand side of the slit in $S_{i+1}$ and the left-hand side of the slit in $S_1$ is glued on the right-hand side of the one in $S_n$.\\
 Denote by $T$ the resulting flat surface together with the canonical projection $\pi:T \to S$ which is a flat $n$-sheeted covering. The branch points are the endpoints of the slits and the lifts of all slits in $T$ to the flat universal cover $\tilde{T}$ is a countable disjoint union of isometrically embedded $n$-valent trees of edge length $l$. Each vertex is the preimage of one of the branch points. Let $B(R)$ be a ball of radius $R$ in $\tilde{T}$ whose center is a vertex. $B(R)$ contains at least $n^{R/l-1}$ vertices.\\ 
One deduces that 
$$e(T)\geq \frac{\log(n)}{l}+\log(1/2)$$ 
Therefore, one cannot expect a growth rate which is smaller than an expression inverse proportional to distance between the branch points and logarithmic in the combinatorics of the cover. We will show that this inequality is almost sharp. Denote by $l_b(T)$ the minimal distance between branch point and by $\lambda(T)$ a measure for the combinatorics i.e. the sum of the number of sheets and the number of branch points. 
\begin{Thm}\label{thmrcov}]
There is some constant $C>0$ and a function $a(S)>0$ that for each flat branched covering $\pi:T \to S$ the
 volume entropy $e(T)$ is bounded by the inequality 
 $$e(T)\leq (e(S)+1) \left(a(S)+\frac{C \cdot \log(\lambda(T))} {l_b(T)}\right)$$ 
The same holds for the Hausdorff dimension of the Gromov boundary.
 \end{Thm}
 The methods to prove the statement are mainly combinatorial. We first show the claim for the special case that the branch points in $T$ project to singularities in $S$. Afterwards we show that the general case can be reduced to the first.\\
We notice the following observations:
\begin{Lem}\label{Lemcomlifts}
 Let $\pi:T \to S$ be a flat branched covering. Denote by $B_T\subset T$ the branch points and $B_S:=\pi(B_T)$ its image. Assume that the branching index is bounded from above by some constant 
$$n:=\sup_{y\in B_T} ind(y)+1$$ 
\begin{enumerate}
 \item For an arc $c:[0,t] \to S $ in $S$ which passes $k$ times through points in $B_S$ and for each preimage $y\in \pi^{-1}(c(0))$ of the starting point of $c(0)$ in $T$ there are at most $k^n$ connected arcs $c'$ in $T$ which project to $c$ and have the same starting point $y$.
\item For an arc $c_0':[0,t] \to T $ in $T$ which passes $k$ times through branch points there are at most $k^n$ connected arcs $c'$ in $T$ which project to $\pi(c_0')$, have the starting point $c'_0(0)$ and pass through $k$ branch points.
\end{enumerate}
\end{Lem}
Recall that by Proposition \ref{Propconsc} for any closed flat surface $S$ there exists a constant $C_l(S)>0$ so that for any two saddle connections $s_1,s_2$ there is a geodesic $g$ of length $l(g) \leq l(s_1)+l(s_2)+C_l(S)$ which first coincides with $s_1$ and eventually with $s_2$.
\begin{Prop}
Let $\pi:T\to S $ be an $n$-sheeted flat covering which branches at most over singularities in $S$ and denote by $N_*(R),~*=S,T$ the counting function on the flat universal cover.\\
Then there is a constant $a_1(S)>0$ so that for $R\geq 1$ it follows
$$N_T(R) \leq (2n)^{R\cdot a_1(S)} N_{S}\left(R\cdot a_1(S)\right)$$
\end{Prop}
\begin{proof}
Denote by $l_0(S)=L_0(T)$ the length of the shortest saddle connection on the space $S,T$ and choose $a_1(S)\geq \frac{2}{l_0(S)}$.\\ 
Fix singularities $x_S\in S,~x_T\in T$ such that $\pi(x_T)=x_S$ and denote by $\mathcal{L}_*(R),~*=S,T $ the set of all locally geodesic parametrized loops of length at most $R$ with the starting point $x_*$. Let $\mathcal{P}_S(R)$ be the set of parametrized loops $h$ in $S$ which have the following properties:
\begin{enumerate}
\item $h$ emanates from $x_S$ and ends at $x_S$.
 \item $h$ is of length of $h$ is at most $R$.
\item $h$ is locally geodesic outside $\Sigma_{S}$.
\end{enumerate}
The projection $\pi: T \to S$ induces a map $\Phi:\mathcal{L}_T(R) \to \mathcal{P}_S(R) $. For a loop $h\in \mathcal{P}_S(R)$ we have to estimate the maximal number of possible preimages.\\
 Recall that the projection $\pi: T\to S$ is an $n$-sheeted branched covering, so the branching index of each branch point is at most $n-1$. The loop $h$ is, away from the singularities, a local geodesic and the starting and endpoint of $h$ is a singularity. Therefore, $h$ is a concatenation of at most $R \cdot a_1(S)$ saddle connections.\\ 
By Lemma \ref{Lemcomlifts} there are at most $n^{R \cdot a_1(S)}$ different geodesic loops in $T$ which emanate from $x_T$ and project to the same loop $h$. So, 
$$N_{T}(R)\leq n^{R\cdot a_1(S)}|\mathcal{P}_S(R)|$$
It remains to compare the cardinality of the set $\mathcal{P}_S(R)$ with the counting function $N_{S}(R)=|\mathcal{L}_S(R)|$. \\
 We define an injective mapping $\Psi$ which maps a piecewise geodesic loop $h \in \mathcal{P}(R)$ to a geodesic loop $ g \in \mathcal{L}\left( R\left(1+\frac{C_l(S)}{l_0(S)}\right)\right)$ which we equip with a combinatorial datum which consists in a coloring of each saddle connection in $g$ with the color red or green.\\ 
As $h\in \mathcal{P}(R)$ is locally geodesic outside the singularities and the endpoints of $h$ are singularities, $h$ is a concatenation of saddle connections $s_1* \ldots *s_m$. Let $s_{i-1}*s_i$ be the incoming and outgoing saddle connections at the singularity $h(t_i)$. By Proposition \ref{Propconsc} there exists a local geodesic $\alpha_i$ whose length is at most $C_l(S)+l(s_{i-1})+l(s_i)$ and which first leaves $s_{i-1}$ and eventually passes through $s_{i}$. If $h$ is not locally geodesic at $h(t_i)$, we replace $s_{i-1}*s_i$ by the local geodesic $\alpha_i$. \\
We can do this construction successively at all points $h(t_i)$ where $h$ is not locally geodesic.
Let $g$ be the resulting locally geodesic loop. The starting point and endpoint of $g$ equals the starting and endpoint of $h$ which is $x_S$.\\ 
 The length of $g$ is bounded from above by:
$$l(g)\leq l(h)+ \frac{RC_l(S)}{l_0(S)}\leq R\left(1 + \frac{C_l(S)}{l_0(S)}\right)$$
The map $h\to g$ is not necessarily injective. It is possible that the same loop $g$ arises from different points $h,h' \in \mathcal{P}(R)$. That is why we add a combinatorial datum. \\
 $h$ is a concatenation of saddle connections $s_1* \ldots *s_m$. $g$ can also be written as such a concatenation $s'_1* \ldots* s'_n$ which arises from $s_1* \ldots *s_m$ by gluing in additional saddle connections. We color each saddle connection $s'_i$ with one of the colors red or green, so that the subsequence of green saddle connections equals the original sequence $s_1*\ldots *s_m$. 
Let $g'$ be the colored loop. The mapping 
$$\Psi(h):=g'$$
is injective, as we can reconstruct $h$ from $g'$ by removing the red saddle connections. 
Since each locally geodesic loop $g$ consists of at most $\frac{l(g)}{l_0(S)}$ saddle connections, $g$ can be colored in at most $2^{\frac{l(g)}{l_0(S)}}$ different ways.\\
 Increase $a_1(S)$ so that $a_1(S)\geq \max\left\{\frac{1}{l_0(S)}\left(1+\frac{C_l(S)}{l_0(S)}\right), 1+\frac{C_l(S)}{l_0(S)}\right\} $ and observe that
$$|\mathcal{P}_S(R)| \leq 2^{R \cdot a_1(S)} \left|\mathcal{L}_S\left( R \cdot a_1(S) \right) \right|$$
The two formulae allow to compare the counting functions for $S$ and $T$.
$$N_T(R) \leq n^{R\cdot a_1(S)} \left(2^{R\cdot a_1(S)}N_{S}\left(R\cdot a_1(S)\right)\right)$$
 \end{proof}
\begin{Cor} \label{CorEntrspeccover}
Let $\pi:T\to S$ be an $n$-sheeted flat branched covering which branches at most over the singularities and let $a_1(S)>0$ be as above. The entropy of $T$ can be estimated in terms of the following formula: 
 $$e(T)\leq a_1(S)(\log(2n) +e(S))$$ 
\end{Cor}
To compute the entropy of general branched coverings, one makes use of the following observation
\begin{Lem} \label{Lemcovergeod}
 Let $\pi:T\to S$ be a flat $n$-sheeted branched covering with exactly one ramification point $x_T$. Assume that $x_T$ ramifies with maximal index $n\geq 3$ over $x_S=\pi(x_T)$ and let $c:[0,t] \to S$ be a connected arc which is locally geodesic outside $x_S$. 
Then, for each $y\in \pi^{-1}(c(0))$ there exists some local geodesic $g$ in $T$ such that $\pi\circ g=c$ and $g(0)=y$.
 \end{Lem}
\begin{Prop}
Let $\pi: T\to S$ be an $n$-sheeted flat branched covering and let $B_T\subset T$ be the set of all branch points in $T$. Denote by $l_b(T)$ the minimal distance between any two branch points on $T$ and let $k:=|B_T|$ be the cardinality of $B_T$.\\
 Fix a singularity $\varsigma_S \in S$ and fix a 3-sheeted flat covering $\pi: T_0 \to S$ which ramifies with maximal index over $\varsigma_S$ and has no other branch points. Then 
$$N_T(R)\leq n^{R/l_b(T)}N_{T_0}\left(R\left(1+\frac{2\mathit{diam}(S)}{l_b(T)}\right)\right)k^{\frac{R}{l_b(T)}}$$
\end{Prop}
\begin{proof}
Fix a point $\varsigma_T \in \pi^{-1}(\varsigma_S) \subset T$ and $\varsigma_{T_0} = \pi^{-1}(\varsigma_S) \in T_0$ the unique preimage of $\varsigma_S$ in $T_0$. Let $B_T=\{b_1\ldots b_k\} \subset T$ be the set of all branch points in $T$ and for each $b_i$ fix a shortest geodesic $c_i$ in $S$ connecting $\pi(b_i)\in S$ with $\varsigma_S$.\\
Let $\mathcal{L}_*(R),~*=T,T_0 $ be the set of parametrized locally geodesic loops of length at most $R$ in $*$, which connect $\varsigma_{*}$ with itself. We define a mapping 
$$\Phi:\mathcal{L}_T(R) \to \mathcal{L}_{T_0}\left(R\left(1+\frac{2\mathit{diam}(S)}{l_b(T)}\right)\right)\times B_T^{\frac{R}{l_b(T)}}$$ 
For a loop $g \in \mathcal{L}_T(R)$ cut $g$ at each branch point $a_i \in B_T$ and obtain local geodesics $g_i, i=1\ldots m,~m\leq \frac{R}{l_b(T)}$ whose projection $\pi\circ g_i\subset S$ is a local geodesic which starts at the image of some branch point $b_{j_i}$ and ends at some other $b_{j_{i+1}}$. Choose geodesics $c_{j_i},c_{j_{i+1}}$ so that $h_i:=c_{j_i}*\pi\circ g_i*c_{j_{i+1}}^{-1}$ is a connected arc with the starting and endpoint $\varsigma_ S$. If $i=1$ resp. $i=m$, one only has to concatenate $c_j$ at the endpoint resp. at the starting point of $\pi \circ g_i$.\\ 
Let $h_{i,\fl} \in [h_i] \subset S$ be the length-minimizing local geodesic in the homotopy class of arcs with fixed endpoints and denote by 
$$h:=h_{1,\fl}*h_{2,\fl}\ldots* h_{m,\fl}$$ 
the concatenation of all such arcs. $h$ is a loop which is locally geodesic outside $\varsigma_S$ and of length 
$$l(h)\leq \sum_{i}(l(g_i)+2l(c_{j_i}))\leq R+\frac{2 R}{l_b(T)}\mathit{diam}(S) $$
By Lemma \ref{Lemcovergeod} we can choose a lift $\alpha$ of $h$ in $T_0$ which is a locally geodesic loop as $\varsigma_{S}$ has exactly one preimage in $T_0$. We define the mapping 
$$\Phi: g \mapsto (\alpha,a_i)$$
where $a_i\in B_T$ is the endpoint of $g_i$.
We have to count the number of possible preimages.\\
Let $g\not=g'\in \mathcal{L}_T(R)$ be different geodesics loops in $T$ that $\Phi(g)=(\alpha,a_i)= \Phi(g')$ and so the data $a_i$ are equal.\\ 
If the projections $\pi(g)\not=\pi(g') \subset S$ are different, by uniqueness of arcs in homotopy classes with fixed endpoint, at least one of the homotopy classes $[g_i]$ and $[g'_i]$ differ and so, the concatenations $h$ and $h'$ are different. This is a contradiction as $\alpha$ cannot be a lift of $h$ and $h'$ at the same time.\\ 
So, $\Phi(g)=\Phi(g')$ only if $\pi(g)=\pi(g')$ and only if the branch point data $a_i$ equal as well.\\
By Lemma \ref{Lemcomlifts} for each $m>0$ there are at most $n^{m}$ different loops $g\subset T$ with the same starting point which pass through $m$ branch points and which project to the same arc $\pi(g)\subset S$.\\
Therefore, the map $\Phi$, restricted to loops of length at most $R$, is at most $n^{R/l_b(T)}$-to-$1$.\\ 
Consequently,
$$N_T(R)\leq n^{\frac{R}{l_b(T)}}N_{T_0}\left(R\left(1+\frac{2\mathit{diam}(S)}{l_b(T)}\right)\right)|B_T|^{\frac{R}{l_b(T)}}$$
\end{proof}
\begin{Cor} \label{Corgencovtospeccover}
Let $\pi:T \to S$ be a flat branched $n$-sheeted covering with $k$ branch point and let $l_b(T)$ be the minimal distance of branch points in $T$. Let $\pi:T_0 \to S$ be defined as above. Then 
$$e(T)\leq \left(1+\frac{2\mathit{\mathit{diam}}(S)}{l_b(T)}\right)e(T_0)+\frac{\log(n)+\log(k)}{l_b(T)}$$
\end{Cor}
Theorem \ref{thmrcov} follows from Corollary \ref{Corgencovtospeccover} and Corollary \ref{CorEntrspeccover}. 

\bibliographystyle{alpha}	
 \bibliography{datab}
\noindent
MATHEMATISCHES INSTITUT DER UNIVERSIT\"AT BONN\\
ENDENICHER ALLEE 60,\\
53115 BONN, GERMANY\\
e-mail: klaus@math.uni-bonn.de

\end{document}